\documentclass[reqno]{amsart}
\usepackage{amsmath}
\numberwithin{equation}{section}
\usepackage{amsfonts}
\usepackage{booktabs}
\usepackage{hyperref}
\hypersetup{
    colorlinks,%
    citecolor=blue,%
    filecolor=blue,%
    linkcolor=blue,%
    urlcolor=blue
}
\usepackage[authoryear]{natbib}
\newtheorem{proposition}{Proposition}[section]
\newtheorem{lemma}[proposition]{Lemma}
\newtheorem{corollary}[proposition]{Corollary}
\newtheorem{theorem}[proposition]{Theorem}

\newcommand{\thlabel}[1]{\label{th:#1}}
\newcommand{\thref}[1]{Theorem~\ref{th:#1}}
\newcommand{\selabel}[1]{\label{se:#1}}
\newcommand{\seref}[1]{Section~\ref{se:#1}}
\newcommand{\lelabel}[1]{\label{le:#1}}
\newcommand{\leref}[1]{Lemma~\ref{le:#1}}

\newcommand{\colabel}[1]{\label{co:#1}}
\newcommand{\coref}[1]{Corollary~\ref{co:#1}}

\newcommand{\eqlabel}[1]{\label{eq:#1}}
\newcommand{\equref}[1]{(\ref{eq:#1})}

\def\et{\eta}

\def\m{{\mu}}
\def\n{{\nu}}

\def\l{{\lambda}}

\def\a{{\alpha}}
\def\e{{\varepsilon}}

\def\G{{\Gamma}}
\def\d{{\delta}}

\def\t{{\theta}}
\def\s{{\sigma}}
\def\z{{\zeta}}
\def\qed{\hfill\Box}

\def\1{{\mathbf 1}}

\def\PP{{\mathbb P}}

\def\MO{{\mathcal O}}

\begin{document}
\bibliographystyle{plainnat}

\title[Distribution of extremes of 1 dependent r.v. sequences]{Approximation for the distribution of extremes of one dependent stationary sequences of random variables}
\author{Am\u arioarei Alexandru}
\address{Laboratoire de Math\'emathiques Paul Painlev\'e, UMR 8524, Univerit\'e de Sciences et Technologies de Lille 1, France}
\address{INRIA Nord Europe/Modal, France}
\address{National Institute of R\&D for Biological Sciences, Bucharest, Romania}
\email{alexandru.amarioarei@inria.fr}
\subjclass[2000]{62E17}
\keywords{extremes, 1-dependent stationary sequences, scan statistics}
\begin{abstract}
{In this  paper we improve some existing results concerning the approximation of the distribution of extremes of a 1-dependent and stationary sequence of random variables. We enlarge the range of applicability and improve the approximation error. An application to the study of the distribution of scan statistics generated by Bernoulli trials is given.}
\end{abstract}

\maketitle
\section{Introduction}\selabel{sec1}
\noindent
The starting point of this paper is a series of results of \cite{Haiman}, concerning the extreme of a 1-dependent and stationary sequence of random variables.
Let $(X_n)_{n\geq1}$ be a sequence of strictly stationary 1-dependent random variables (for any $t\geq1$ we have $\s(X_1,\dots,X_t)$ and $\s(X_{t+2},\dots)$ are independent) with marginal distribution function $F(x)=\PP(X_1\leq x)$. Let $x$ such that
$$
\inf\{u|F(u)>0\}<x<\sup\{u|F(u)<1\}.
$$
Define the sequences
\begin{eqnarray}
p_n&=&p_n(x)=\PP(\min\{X_1,X_2,\dots,X_n\}>x),\ \ n\geq1,\ p_0=1, \eqlabel{eq1.1}\\
q_n&=&q_n(x)=\PP(\max\{X_1,X_2,\dots,X_n\}\leq x),\ \ n\geq1 \eqlabel{eq1.2}
\end{eqnarray}
and the series
\begin{equation}
C(z)=C_x(z)=1+\displaystyle\sum_{k=1}^{\infty}{(-1)^kp_{k-1}z^k}.\eqlabel{eq1.3}
\end{equation}
In \cite{Haiman}, the author proved the following results:
\begin{theorem}\thlabel{th1}
For $x$ such that $0<p_1(x)\leq0.025$, $C_x(z)$ has a unique zero $\l(x)$, of order of multiplicity $1$, inside the interval $(1,1+2p_1)$, such that
\begin{equation}
|\l-(1+p_1-p_2+p_3-p_4+2p_1^2+3p_2^2-5p_1p_2)|\leq87p_1^3. \eqlabel{eqth1}
\end{equation}
\end{theorem}
\noindent
\begin{theorem}\thlabel{th2}
We have
$$
q_1=1-p_1,\ \ \ \ q_2=1-2p_1+p_2,\ \ \ \ q_3=1-3p_1+2p_2+p_1^2-p_3
$$
and for $n>3$ if $p_1\leq0.025$,
\begin{equation}
|q_n\l^n-(1-p_2+2p_3-3p_4+p_1^2+6p_2^2-6p_1p_2)|\leq 561p_1^3.\eqlabel{eqth2}
\end{equation}
\end{theorem}
\noindent
These results were successfully applied in a series of applications: the distribution of the maximum of the increments of the Wiener process (\cite{Haiman}), extremes of Markov sequences (\cite{Haiman1}), the distribution of scan statistics, both in one dimensional (see \cite{Haiman2,Haiman2007}) and two dimensional case (see \cite{HaimanPreda1,HaimanPreda2}).\\ \noindent
Following the same lines of proofs as in \cite{Haiman}, we improve \thref{th1} and \thref{th2} by enlarging the range of applicability and providing sharper error bounds. The main results are presented in \seref{sec2}. In \seref{sec3}, we present an application to the study of the distribution of one dimensional discrete scan statistics emphasizing the difference between the new and the old results. Proofs are presented in \seref{sec4}.
\section{Main results}\selabel{sec2}
\noindent
The following theorem gives a parametric form of the \thref{th1} improving both the range of $p_1(x)$, from $0.025$ to $0.1$, and the error coefficient:
\begin{theorem}\thlabel{T1}
For $x$ such that $0<p_1(x)\leq\a\leq0.1$, $C_x(z)$ has an unique zero $\l(x)$, of order of multiplicity $1$, inside an interval of the form $(1,1+lp_1)$, such that
\begin{equation}\eqlabel{eqT1.1}
|\l-(1+p_1-p_2+p_3-p_4+2p_1^2+3p_2^2-5p_1p_2)|\leq K(\a)p_1^3,
\end{equation}
where $l=l(\a)>t_2^3(\a)$, $t_2(\a)$ is the second root in magnitude of the equation $\a t^3-t+1=0$ and $K(\a)$ is given by
\begin{equation}\eqlabel{eqT1.2}
K(\a)=\frac{\frac{11-3\a}{(1-\a)^2}+2l(1+3\a)\frac{2+3l\a-\a(2-l\a)(1+l\a)^2}{\left[1-\a(1+l\a)^2\right]^3}}{1-\frac{2\a(1+l\a)}{\left[1-\a(1+l\a)^2\right]^2}}.
\end{equation}
\end{theorem}
\noindent
Using the properties of the 1-dependent sequence, an immediate consequence of \thref{T1} is the following:
\begin{corollary}\colabel{C1}
Let $\l$ be defined as in \thref{T1}, then
\begin{equation}\eqlabel{eqC1}
|\l-(1+p_1-p_2+2(p_1-p_2)^2)|\leq (1+\a K(\a)) p_1^2.
\end{equation}
\end{corollary}
\noindent
To get a better grasp of the bounds in \thref{T1} and \coref{C1}, we present, for selected values of $\a$, the values taken by the coefficients in Eq.\equref{eqT1.1} and Eq.\equref{eqC1}:
\begin{table}[ht]
\centering
 \begin{tabular}{rrrr}
\toprule
$\a$  & $l$  & $K(\a)$  & $1+\a K(\a)$  \\
\midrule
$0.100$ & $1.5347$  &  $38.6302 $ & $4.8630 $   \\


$0.050$ & $1.1893$  &  $21.2853$ & $ 2.0642$  \\

$\mathbf{0.025}$ & $\mathbf{1.0835}$  &  $\mathbf{17.5663}$ & $\mathbf{1.4391}$ \\

$0.010$ & $1.0313$  &  $15.9265$ & $1.1592$  \\


\bottomrule
\end{tabular}
\caption{Selected values for the error coefficients in \thref{T1} and \coref{C1}}
\label{tab1}
\end{table}

\noindent
Notice that for the value considered in \cite{Haiman}, i.e. $\a=0.025$, our corresponding value for the error coefficient in Eq.\equref{eqT1.1} is almost five times smaller than in Eq.\equref{eqth1}. The following result improves \thref{th2}:
\noindent

\begin{theorem}\thlabel{T2}
Lets suppose that $x$ is such that $0<p_1(x)\leq\a\leq0.1$ and define $\et=1+l\a$ with $l=l(\a)>t_2^3(\a)$ and $t_2(\a)$ the second root in magnitude of the equation $\a t^3-t+1=0$. If $\l=\l(x)$ is the zero obtained in \thref{T1}, then the following relation holds
\begin{equation}\eqlabel{eqT2.1}
|q_n\l^n-(1-p_2+2p_3-3p_4+p_1^2+6p_2^2-6p_1p_2)|\leq \G(\a)p_1^3,
\end{equation}
where $\G(\a)=L(\a)+E(\a)$, $K(\a)$ is given by Eq.\equref{eqT1.2} and
\begin{eqnarray}
L(\a)&=&3K(\a)(1+\a+3\a^2)[1+\a+3\a^2+K(\a)\a^3]+\a^6K^3(\a)\nonumber\\
     & &+9\a(4+3\a+3\a^2)+55 \eqlabel{eqP1}\\
E(\a)&=&0.1+\frac{\et^5\left[1+(1-2\a)\et\right]^4\left[1+\a(\et-2)\right]\left[1+\et+(1-3\a)\et^2\right]}{2(1-\a\et^2)^4\left[(1-\a\et^2)^2-\a\et^2(1+\et-2\a\et)^2\right]}\eqlabel{eqE1}
\end{eqnarray}
\end{theorem}
\noindent
The next corollary is an immediate consequence of \thref{T2}:
\begin{corollary}\colabel{C2}
In the conditions of \thref{T2} we have
\begin{equation}
|q_n\l^n-(1-p_2)|\leq(3+\a\G(\a))p_1^2. \eqlabel{eqC2}
\end{equation}
\end{corollary}

\noindent
The error coefficient in Eq.\equref{eqT2.1} is smaller in comparison with the corresponding one from Eq.\equref{eqth2}, as the following table can show:
\begin{table}[h]
\centering
 \begin{tabular}{rrr}
 \toprule
$\a$  &$\G(\a)$  & $3+\a \G(\a)$  \\
\midrule
$0.100$ &  $480.696 $ & $51.0696 $   \\


$0.050$ &  $180.532$ & $ 12.0266$  \\

$\mathbf{0.025}$ &  $\mathbf{145.202}$ & $\mathbf{6.6300}$ \\

$0.010$ &   $131.438$ & $4.3143$  \\


\bottomrule
\end{tabular}
\caption{Selected values for the error coefficients in \thref{T2} and \coref{C2}}
\label{tab2}
\end{table}

\noindent
Combining the results obtained in \thref{T1} and \thref{T2} we get the following approximation:
\begin{theorem}\thlabel{T3}
Let $x$ such that $q_1(x)\geq1-\a\geq 0.9$. If $\G(\a)$ and $K(\a)$ are the same as in \thref{T2}, then
\begin{equation}
\left|q_n-\frac{6(q_1-q_2)^2+4q_3-3q_4}{(1+q_1-q_2+q_3-q_4+2q_1^2+3q_2^2-5q_1q_2)^n}\right|\leq\Delta_1(1-q_1)^3 \eqlabel{eqT3}
\end{equation}
with
\begin{equation}
\Delta_1=\Delta_1(\a,n)=\G(\a)+nK(\a). \eqlabel{eqT3.2}
\end{equation}
\end{theorem}

\noindent
In the same fashion combining the results from \coref{C1} and \coref{C2} we get
\begin{theorem}\thlabel{T4}
If $x$ is such that $q_1(x)\geq1-\a\geq 0.9$, then
\begin{equation}
\left|q_n-\frac{2q_1-q_2}{\left[1+q_1-q_2+2(q_1-q_2)^2\right]^n}\right|\leq\Delta_2(1-q_1)^2 \eqlabel{eqT4}
\end{equation}
with
\begin{equation}
\Delta_2=\Delta_2(\a,n,q_1)=3+\G(\a)(1-q_1)+n\left[1+K(\a)(1-q_1)\right], \eqlabel{eqT4.2}
\end{equation}
and where $\G(\a)$ and $K(\a)$ are the same as in \thref{T2}.
\end{theorem}
\noindent
\section{Application to the distribution of scan statistics}\selabel{sec3}
\noindent
In many applications the decision makers have to determine if a certain accumulation of events is \textit{normal} or not, where by \textit{normal} we mean that it can be explained by an underlying probability model defined by a null hypothesis of randomness. One way of determining if such a cluster of observations is exceptionally or not is by using the scan statistics.
For a comprehensive methodological treatment and a rich source of applications of this test statistic, one may study the books of \cite{Glaz2001} and more recently the one of \cite{Glaz2009}.\\
\noindent
Let $Y_1, Y_2, \dots, Y_N$ be a sequence of independent and identically distributed (i.i.d.) random variables and $m$, $1\leq m\leq N$, be a fixed positive integer. If we consider the random variables
\begin{equation}\eqlabel{seq1}
Z_t=\displaystyle\sum_{i=t}^{t+m-1}{Y_i},\ \ 1\leq t\leq N-m+1
\end{equation}
then the one dimensional discrete scan statistic is defined by
\begin{equation}\eqlabel{seq2}
S_m(N)=\displaystyle\max_{1\leq t\leq N-m+1}{Z_t}.
\end{equation}
To get an intuitive meaning of the above definition, if $Y_i$ are integer valued random variables then we can interpret them as the number of observed events at the time $i$. The scan statistics is then viewed as the maximum number of events observed in any contiguous period of length $m$ within the interval $\{1,2,\dots,N\}$. Since exact formulas for the distribution of $S_m$ exist only in a small number of situations (see for example \cite{Glaz2001}, Chapter 13), various approximation methods and bounds have been proposed. In what follows we will use the approximation method developed in \cite{Haiman2007} but with the help of the results obtained in \seref{sec2}. The method is based on the important observation that, in the i.i.d. case, the discrete scan statistic can be expressed as an extreme of a 1-dependent stationary sequence. It is easy to see that the random variables
\begin{equation}\eqlabel{seq3}
W_k=\displaystyle\max_{(k-1)m+1\leq s\leq km+1}{Z_s},\ \ \ k=1,2,\dots
\end{equation}
form a 1-dependent stationary sequence and that the following relation holds,
\begin{equation}\eqlabel{seq4}
S_m\left(Lm\right)=\max{(W_1,W_2,\dots,W_{L-1})}.
\end{equation}
Observe that for general $N$, one can consider $L=[N/m]$ and then apply the inequality
\begin{equation}\eqlabel{seq5}
\PP\left(S_m[(L+1)m]\leq n\right)\leq \PP\left(S_m(N)\leq n\right)\leq \PP\left(S_m(Lm)\leq n\right).
\end{equation}
When $\PP(W_1>n)\leq 0.1$, we can apply the results from \thref{T4} to the sequence $W_1,W_2,\dots,W_{L-1}$ to obtain the following approximation for the distribution of the scan statistic $S_m$:
\begin{equation}\eqlabel{seq6}
\PP\left(S_m(Lm)\leq n\right)\approx \frac{2q_1-q_2}{\left[1+q_1-q_2+2(q_1-q_2)^2\right]^{(L-1)}},
\end{equation}
with an error bound of about
\begin{equation}\eqlabel{seq7}
E=\left\{3+\G(\a)(1-q_1)+(L-1)\left[1+K(\a)(1-q_1)\right]\right\}(1-q_1)^2,
\end{equation}
where $\a=\PP(W_1>n)$ and $q_1,q_2$ are defined by Eq.\equref{eq1.2} as
\begin{eqnarray}
q_1&=&\PP\left(W_1\leq n\right)=\PP(S_m(2m)\leq n), \eqlabel{seq8}\\
q_2&=&\PP\left(W_1\leq n,W_2\leq n\right)=\PP(S_m(3m)\leq n). \eqlabel{seq9}
\end{eqnarray}
We should mention that in \citet[Theorem 4]{Haiman}, the author obtained the following formula for the approximation error
\begin{equation}\eqlabel{seq10}
EH=\left\{9+561(1-q_1)+3.3(L-1)\left[1+4.7(L-1)(1-q_1)^2\right]\right\}(1-q_1)^2.
\end{equation}
Notice that if $q_1,q_2,q_3,q_4$ are known then we can apply \thref{T3} to get a better approximation for the distribution of $S_m$, even though in most applications Eq.\equref{seq6} will suffice.\\
Next we will restrict ourselves to the case of Bernoulli $0-1$ process, that is $Y_i$'s are Bernoulli trials with $\PP(Y_i=1)=p=1-\PP(Y_i=0)$. In this particular framework, Naus (see for example \cite{Glaz2001}, Theorem 13.1) provided exact formulas for $\PP(S_m(N)\leq n)$ when $N=2m$ and $N=3m$, namely for $q_1$ and $q_2$ from Eq.\equref{seq8} and Eq.\equref{seq9}. The following tables illustrates a compared study between the error formula given by Eq.\equref{seq7} and the corresponding error bound used by \cite{Haiman2007} along with the exact and approximated value for the distribution of the scan statistic:
\begin{table}[ht]
\centering
 \begin{tabular}{cccccccc}
\toprule
$n$  & $q_1$  & $q_2$ & $Approx$ & $Exact$ & $EH$ & $E$   \\
  & {\scriptsize Eq.\equref{seq8} } &{\scriptsize Eq.\equref{seq9}} & {\scriptsize Eq.\equref{seq6}} &  & {\scriptsize Eq.\equref{seq10} }& {\scriptsize Eq.\equref{seq7} } \\
\midrule
$2$ & $0.97131$  &  $0.95181$ & $0.82715$ & $0.82582$ & $-$ & $0.01712$    \\
$3$ & $0.99716$  &  $0.99500$ & $0.98001$ & $0.98000$ & $0.00032$ & $0.00010$    \\
$4$ & $0.99982$  &  $0.99967$ & $0.99865$ & $0.99865$ & $1\times10^{-6}$ & $3\times10^{-7}$    \\
$5$ & $0.99999$  &  $0.99998$ & $0.99994$ & $0.99994$ & $2\times10^{-9}$ & $6\times10^{-10}$    \\
$6$ & $1.$  &  $1.$ & $0.99999$ & $0.99999$ & $1\times10^{-12}$ & $4\times10^{-13}$    \\
$7$ & $1.$  &  $1.$ & $1.$ & $1.$ & $3\times10^{-16}$ & $9\times10^{-17}$    \\
\bottomrule
\end{tabular}
\caption{Distribution of the scan statistic $\PP\left(S_m(Lm)\leq n\right)$ for $m=9$, $p=0.05$, $L=10$}
\label{tab3}
\end{table}
\par
\begin{table}[ht]
\centering
 \begin{tabular}{ccccccc}
\toprule
$n$  & $q_1$  & $q_2$ & $Approx$ & $Exact$ & $EH$ & $E$   \\
  &{\scriptsize Eq.\equref{seq8} } &{\scriptsize Eq.\equref{seq9}} & {\scriptsize Eq.\equref{seq6}} &  & {\scriptsize Eq.\equref{seq10}} & {\scriptsize Eq.\equref{seq7}  } \\
\midrule
$1$ & $0.96860$  &  $0.94910$ & $0.74617$ & $0.74353$ & $-$ & $0.02927$    \\
$2$ & $0.99813$  &  $0.99677$ & $0.98061$ & $0.98060$ & $0.00019$ & $0.00006$    \\
$3$ & $0.99993$  &  $0.99987$ & $0.99922$ & $0.99922$ & $2\times10^{-7}$ & $8\times10^{-8}$    \\
$4$ & $0.99999$  &  $0.99999$ & $0.99998$ & $0.99998$ & $1\times10^{-10}$ & $4\times10^{-11}$    \\
$5$ & $1.$  &  $1.$ & $1.$ & $1.$ & $4\times10^{-14}$ & $1\times10^{-14}$    \\
\bottomrule
\end{tabular}
\caption{Distribution of the scan statistic $\PP\left(S_m(Lm)\leq n\right)$ for $m=10$, $p=0.0165$, $L=15$}
\label{tab4}
\end{table}
\par
\noindent
The exact values for the distribution of the scan statistics presented in the tables ~\ref{tab3} and ~\ref{tab4} (column "Exact") are computed using the Markov chain embedding technique described in \cite{Fu}.
\section{Proofs of the results}\selabel{sec4}
\noindent
\subsection{Proof of \thref{T1}}
\noindent
The proof will follow closely that of \cite{Haiman}. Using the stationarity and 1-dependence of the sequence we have
\begin{eqnarray*}
p_n=\PP(X_1>x,X_2>x,\dots,X_n>x)&\leq&\PP(X_1>x,X_3>x,\dots,X_n>x)\\
                                &=&\PP(X_1>x)\PP(X_3>x,\dots,X_n>x)\\
                                &=&p_1p_{n-2},
\end{eqnarray*}
which give the basic inequality
\begin{equation}\eqlabel{t1eq1}
p_n\leq p_1^{[\frac{n+1}{2}]}.
\end{equation}
To show that $C(z)$ has a zero in the interval $(1,1+lp_1)$ it is enough to prove that $C(1)>0$ and $C(1+lp_1)<0$. It easy to see that $C(1)>0$, since
\begin{equation}\eqlabel{t1eq2}
C(1)=1+\displaystyle\sum_{k=1}^{\infty}{(-1)^kp_{k-1}}=\underbrace{(p_1-p_2)}_{\geq0}+\underbrace{(p_3-p_4)}_{\geq0}+\dots\geq0
\end{equation}
For $C(1+lp_1)$ we have:
\begin{eqnarray}
C(1+lp_1)&=&1+\displaystyle\sum_{k=1}^{\infty}{(-1)^kp_{k-1}(1+lp_1)^k}\nonumber\\
         &=&-lp_1+\displaystyle\sum_{k=1}^{\infty}{(1+lp_1)^{2k}\underbrace{\left[p_{2k-1}-p_{2k}(1+lp_1)\right]}_{\leq p_{2k-1}\leq p_1^k}}\nonumber\\
         &\leq&-lp_1+\displaystyle\sum_{k=1}^{\infty}{[(1+lp_1)^{2}p_1]^k}.\eqlabel{t1eq3}
\end{eqnarray}
It is easy to see that if $(1+lp_1)^3<l$ then $p_1(1+lp_1)^2<\frac{lp_1}{1+lp_1}<1$, so series in Eq.\equref{t1eq3} is convergent and $C(1+lp_1)<0$. From the definition of $l$ and the relation $1<t_2(\a)<t_2(\a)+\e<\frac{1}{\sqrt{3\a}}$, we obtain ($t_2=t_2(\a)$)
\begin{eqnarray*}
t_2^2+t_2\sqrt[3]{t_2^3+\e}+\sqrt[3]{\left(t_2^3+\e\right)^2}&\leq&t_2^2+t_2(t_2+\e)+\left(t_2+\e\right)^2\\
                                                                           &<&\frac{1}{3\a}+\frac{1}{3\a}+\frac{1}{3\a}=\frac{1}{\a}
\end{eqnarray*}
which imply that $t_2-\sqrt[3]{t_2^3+\e}+\a\e<0$. Combining this last relation with the fact that $t_2$ is a root of the equation $\a t^3-t+1=0$, we have that $(1+lp_1)^3<l$.
For showing that the zero is unique we will prove that $C'(z)$ is strictly decreasing on the interval $(1,1+lp_1)$, i.e. $C'(z)<0$. \noindent
Using Lagrange theorem on $[1,1+a]$ we get for $\t\in(0,1)$
\begin{equation}\eqlabel{t1eq4}
C'(1+a)=C'(1)+aC''(1+\t a)
\end{equation}
We will approximate both $C'(1)$ and $C''(1+\t a)$ as follows:
\begin{equation}\eqlabel{t1eq5}
C'(1)=-1+2p_1-3p_2+R
\end{equation}
where
\begin{equation}\eqlabel{t1eq6}
R=\displaystyle\sum_{k=2}^{\infty}{2k\left(p_{2k-1}-p_{2k}\right)}-\sum_{k=2}^{\infty}{p_{2k}}.
\end{equation}
To approximate $R$ notice that $p_{2k-1}-p_{2k}\geq0$, which implies
\begin{equation*}
-\displaystyle\sum_{k=2}^{\infty}{p_{2k}}\leq R\leq \displaystyle\sum_{k=2}^{\infty}{2k\left(p_{2k-1}-p_{2k}\right)}
\end{equation*}
and using Eq.\equref{t1eq1} we get
\begin{equation}\eqlabel{t1eq7}
-\frac{p_1^2}{1-p_1}=-\displaystyle\sum_{k=2}^{\infty}{p_1^k}\leq R\leq 2p_1\displaystyle\sum_{k=2}^{\infty}{kp_1^{k-1}}=2p_1^2\left[\frac{1}{(1-p_1)^2}+\frac{1}{1-p_1}\right]
\end{equation}
so that
\begin{equation}\eqlabel{t1eq8}
|R|\leq 2p_1^2\left[\frac{1}{(1-p_1)^2}+\frac{1}{1-p_1}\right].
\end{equation}
For $C''(z)$ we have:
\begin{equation}\eqlabel{t1eq9}
C''(z)=2p_1-2\cdot3p_2z+3\cdot4p_3z^2-4\cdot5p_4z^3+5\cdot6p_5z^4-\dots
\end{equation}
Using $z\in(1,1+lp_1)$ and Eq.\equref{t1eq1} we have
\begin{eqnarray}\eqlabel{t1eq10}
C''(z)&\leq& p_1\displaystyle\sum_{k=0}^{\infty}{(2k+1)(2k+2)(z\sqrt{p_1})^{2k}}\nonumber\\
      &\leq& 2p_1\frac{1+3p_1(1+lp_1)^2}{\left[1-p_1(1+lp_1)^2\right]^3}.
\end{eqnarray}
Since $z\leq1+lp_1$, Eq.\equref{t1eq9} and Eq.\equref{t1eq1} implies
\begin{eqnarray}\eqlabel{t1eq11}
C''(z)&\geq&-lp_1\displaystyle\sum_{k=0}^{\infty}{(2k+1)(2k+2)p_{2k+2}z^{2k}}-4\displaystyle\sum_{k=0}^{\infty}{(k+1)p_{2k+2}z^{2k+1}}\nonumber\\
      &\geq&-lp_1^2\displaystyle\sum_{k=0}^{\infty}{(2k+1)(2k+2)(z\sqrt{p_1})^{2k}}-4p_1z\displaystyle\sum_{k=0}^{\infty}{(k+1)(p_1z^2)^k}\nonumber\\
      &\geq&-2lp_1^2\frac{1+3p_1(1+lp_1)^2}{\left[1-p_1(1+lp_1)^2\right]^3}-\frac{4p_1(1+lp_1)}{\left[1-p_1(1+lp_1)^2\right]^2},
\end{eqnarray}
which in relation with Eq.\equref{t1eq10} shows that
\begin{equation}\eqlabel{t1eq12}
|C''(z)|\leq 2lp_1^2\frac{1+3p_1(1+lp_1)^2}{\left[1-p_1(1+lp_1)^2\right]^3}+\frac{4p_1(1+lp_1)}{\left[1-p_1(1+lp_1)^2\right]^2}.
\end{equation}
Combining Eqs.\equref{t1eq4}, \equref{t1eq5}, \equref{t1eq8} and \equref{t1eq11} we can show that $C'(z)<0$ if the following inequality is true
\begin{equation}\eqlabel{t1eq13}
-1+\frac{2p_1}{(1-p_1)^2}+2lp_1^2\frac{1+3p_1(1+lp_1)^2}{\left[1-p_1(1+lp_1)^2\right]^3}<0.
\end{equation}
We observe that the expression on the left hand side of Eq.\equref{t1eq13} is increasing in $p_1$ and since for $p_1=0.1$ we have $l\leq 1.1535$ we get
\begin{equation}\eqlabel{t1eq14}
\frac{2p_1}{(1-p_1)^2}+2lp_1^2\frac{1+3p_1(1+lp_1)^2}{\left[1-p_1(1+lp_1)^2\right]^3}< 0.3
\end{equation}
which verifies that $C'(z)<0$.
\par
\noindent
Now we try to approximate the zero $\l$. From Lagrange theorem applied on the interval $[1,\l]$ we have $C(\l)-C(1)=(\l-1)C'(u)$, with $u\in(1,\l)\subset(1,1+lp_1)$. Since $C(\l)=0$ we get
\begin{equation}\eqlabel{t1eq15}
\l-1=-\frac{C(1)}{C'(u)}
\end{equation}
and taking $\m=p_1-p_2+p_3-p_4+2p_1^2+3p_2^2-5p_1p_2$ as in \cite{Haiman} we obtain
\begin{equation}\eqlabel{t1eq16}
\l-(1+\m)=-\frac{C(1)+\m C'(u)}{C'(u)}.
\end{equation}
Applying Lagrange theorem one more time as in Eq.\equref{t1eq4} and combining with Eq.\equref{t1eq5}, the relation in Eq.\equref{t1eq16} becomes
\begin{equation}\eqlabel{t1eq17}
\l-(1+\m)=-\frac{C(1)+\m(-1+2p_1-3p_2)+\m(R+aC''(1+\t a))}{C'(u)},
\end{equation}
where $a=u-1$ and $\t\in(0,1)$.
\par
\noindent
If we denote $T_1=C(1)+\m(-1+2p_1-3p_2)$, then
\begin{eqnarray}\eqlabel{t1eq18}
T_1&=&(p_1-p_2+p_3-p_4+p_5-p_6+\dots)-\m+(2p_1-3p_2)\m \nonumber\\
    &=&(p_5-p_6+\dots)+(p_1-p_2)(2p_1-3p_2)^2+2(p_1-p_2)(p_3-p_4)\nonumber\\
    & &-p_2(p_3-p_4).
\end{eqnarray}
Using Eq.\equref{t1eq1} in Eq.\equref{t1eq18} we observe that
\begin{equation*}
-p_1^3\leq-p_2(p_3-p_4)\leq T_1\leq\displaystyle\sum_{k=2}^{\infty}{p_1^{k+1}}+4p_1^3+2p_1^3
\end{equation*}
which gives us
\begin{equation}\eqlabel{t1eq19}
|C(1)+\m(-1+2p_1-3p_2)|\leq p_1^3\left[6+\frac{1}{1-p_1}\right].
\end{equation}
Also we can notice that $\m\geq0$ and
\begin{eqnarray}\eqlabel{t1eq20}
\m  &=&(1-p_2)(p_1-p_2)+p_3-p_4+2(p_1-p_2)^2\nonumber\\
    &\leq& (1-p_2)(p_1-p_2)+p_1(p_1-p_2)+2(p_1-p_2)^2\nonumber\\
    &=& 3(p_1-p_2)^2+p_1-p_2\leq p_1(1+3p_1).
\end{eqnarray}
The last step is to find an upper bound for $|C'(u)|^{-1}$. For this observe that
\begin{eqnarray}\eqlabel{t1eq21}
|C'(u)|^{-1}&=&|1-(2p_1u-3p_2u^2+4p_3u^3-5p_4u^4+\dots)|^{-1}\nonumber\\
            &\leq&\left|1-|2p_1u-3p_2u^2+4p_3u^3-5p_4u^4+\dots|\right|^{-1}
\end{eqnarray}
where we used the inequality $|1-x|\geq|1-|x||$. Taking into account that $u\in(1,\l)\subset(1,1+lp_1)$ and denoting the expression inside the second absolute value in the denominator of Eq.\equref{t1eq21} by $T_2$, we have
\begin{eqnarray}\eqlabel{t1eq22}
T_2&=&\displaystyle\sum_{k=1}^{\infty}{2k(p_{2k-1}-p_{2k}u)u^{2k-1}}-\sum_{k=1}^{\infty}{p_{2k}u^{2k}}\nonumber\\
    &\leq&2(p_1-p_2)u\displaystyle\sum_{k=1}^{\infty}{k(p_1u^2)^{k-1}}\leq\frac{2p_1u}{(1-p_1u^2)^2}.
\end{eqnarray}
In the same way
\begin{eqnarray}\eqlabel{t1eq23}
T_2&\geq&-\displaystyle\sum_{k=1}^{\infty}{p_{2k}u^{2k}}-2lp_1\displaystyle\sum_{k=1}^{\infty}{kp_{2k}u^{2k-1}}\nonumber\\
    &\geq&-\displaystyle\sum_{k=1}^{\infty}{p_1^{k}u^{2k}}-2lp_1^2u\displaystyle\sum_{k=1}^{\infty}{k(p_1u^2)^{k-1}}\nonumber\\
    &\geq&-p_1u\frac{2lp_1+u(1-p_1u^2)}{\left(1-p_1u^2\right)^2}>-\frac{2p_1u}{(1-p_1u^2)^2}.
\end{eqnarray}
Combining Eq.\equref{t1eq22} with Eq.\equref{t1eq23} along with $u\leq1+lp_1$, we have
\begin{equation}\eqlabel{t1eq23_2}
|T_2|\leq \frac{2p_1(1+lp_1)}{\left[1-p_1(1+lp_1)^2\right]^2}.
\end{equation}
Substituting Eq.\equref{t1eq23_2} in Eq.\equref{t1eq21} we obtain the bound
\begin{equation}\eqlabel{t1eq24}
|C'(z)|^{-1}\leq\frac{1}{1-\frac{2p_1(1+lp_1)}{\left[1-p_1(1+lp_1)^2\right]^2}}.
\end{equation}
Combining the Eqs.\equref{t1eq8}, \equref{t1eq12}, \equref{t1eq19}, \equref{t1eq20}, \equref{t1eq24} along with the fact that $|a|\leq lp_1$ in Eq.\equref{t1eq17}, we obtain
\begin{eqnarray}\eqlabel{t1eq25}
|\l-(1+\m)|&\leq&\frac{|C(1)+\m(-1+2p_1-3p_2)|+|\m|(|R|+|a||C''(1+\t a)|)}{|C'(u)|}\nonumber\\
            &\leq&K(p_1)p_1^3
\end{eqnarray}
where
\begin{equation}\eqlabel{t1eq26}
K(p_1)=\frac{\frac{11-3p_1}{(1-p_1)^2}+2l(1+3p_1)\frac{2+3lp_1-p_1(2-lp_1)(1+lp_1)^2}{\left[1-p_1(1+lp_1)^2\right]^3}}{1-\frac{2p_1(1+lp_1)}{\left[1-p_1(1+lp_1)^2\right]^2}}.
\end{equation}
To obtain $K(\a)$ it is enough to substitute $p_1$ in the above relation with $\a$ with the additional remark that $l=l(\a)$.
$\qed$

\subsection{Proof of \thref{T2}}
\noindent
For completeness we will give a detailed proof of the theorem even if we repeat most of its ideas from \cite{Haiman}. Remembering that
\begin{equation}\eqlabel{t2eq1}
C(z)=\displaystyle\sum_{k=0}^{\infty}{c_kz^k}=1+\displaystyle\sum_{k=1}^{\infty}{(-1)^kp_{k-1}z^k}
\end{equation}
we define
\begin{equation}\eqlabel{t2eq2}
D(z)=\frac{1}{C(z)}=\displaystyle\sum_{k=0}^{\infty}{d_kz^k}
\end{equation}
which exists since $c_0=1$, and from $C(z)D(z)=1$ we have that $d_0=1$ and
\begin{equation}\eqlabel{t2eq3}
\displaystyle\sum_{j=0}^{n}{d_jc_{n-j}=\displaystyle\sum_{j=0}^{n}{(-1)^{n-j}p_{n-j-1}d_{j}}=0},\ \ n\geq1.
\end{equation}
Now if we define $A_k=\{X_1\leq x,\dots,X_k\leq x,X_{k+1}>x,\dots,X_n>x\}$, for $k=\overline{0,n}$, we obtain $\PP(A_0)=p_n$, $\PP(A_n)=q_n$ and
\begin{equation}\eqlabel{t2eq4}
\PP(A_k)+\PP(A_{k-1})=p_{n-k}q_{k-1}, \ \ k\geq1
\end{equation}
Summing Eq.\equref{t2eq4} over $k$, we deduce that
\begin{equation}\eqlabel{t2eq5}
q_n=\displaystyle\sum_{k=0}^{n}{(-1)^{n-k}p_{n-k}q_{k-1}},\ \ n\geq0,\ \ q_{-1}=q_0=1
\end{equation}
and comparing with Eq.\equref{t2eq3}, after using mathematical induction, we conclude that $d_{n+1}=q_n$ and that
\begin{equation}\eqlabel{t2eq6}
D(z)=\displaystyle\sum_{k=0}^{\infty}{q_{k-1}z^k},\ \ q_{-1}=q_0=1.
\end{equation}
Taking $\l$ as in \thref{T1}, we can write $C(z)=U(z)\left(1-\frac{z}{\l}\right)$ and observe that if we let $U(z)=\displaystyle\sum_{k=0}^{n}{u_kz^k}$ then
\begin{equation}\eqlabel{t2eq7}
\left(\displaystyle\sum_{k=0}^{n}{u_kz^k}\right)\left(1-\frac{z}{\l}\right)=1+\displaystyle\sum_{k=1}^{\infty}{(-1)^kp_{k-1}z^k}
\end{equation}
which shows that
\begin{equation}\eqlabel{t2eq8}
u_n-\frac{u_{n-1}}{\l}=(-1)^np_{n-1},\ u_0=1,\ n\geq1.
\end{equation}
Multiplying Eq.\equref{t2eq8} with $\l^n$ and summing over $n$ we find
\begin{equation}\eqlabel{t2eq9}
u_n=\frac{1+\displaystyle\sum_{k=1}^{n}{(-1)^kp_{k-1}\l^k}}{\l^n},\ \ n\geq1.
\end{equation}
If we denote with $\frac{1}{U(z)}=\displaystyle\sum_{k=0}^{\infty}{t_kz^k}=T(z)$ then
\begin{equation}\eqlabel{t2eq10}
D(z)\left(1-\frac{z}{\l}\right)=T(z)
\end{equation}
and using the same argument as above we get $t_0=d_0=1$ and
\begin{equation}\eqlabel{t2eq11}
t_n=d_n-\frac{d_{n-1}}{\l}, \ \ n\geq1
\end{equation}
so $d_n\l^n=t_0+t_1\l+\dots+t_n\l^n$, that is
\begin{equation}\eqlabel{t2eq12}
q_n\l^{n+1}=t_0+t_1\l+\dots+t_n\l^n+t_{n+1}\l^{n+1}.
\end{equation}
To obtain the desired result we begin by giving an approximation of $u_n$:
\begin{eqnarray}\eqlabel{t2eq13}
|u_n|&=&\frac{\left|1+\displaystyle\sum_{k=1}^{n}{(-1)^kp_{k-1}\l^k}\right|}{\l^n}\overset{C(\l)=0}{=}\frac{\left|\displaystyle\sum_{k=n+1}^{\infty}{(-1)^kp_{k-1}\l^k}\right|}{\l^n}\nonumber\\
     &\leq&\frac{\l^{n+1}}{\l^n}\left|p_n-p_{n+1}\l+p_{n+2}\l^2-p_{n+3}\l^3+\dots\right|\nonumber\\
     &\leq&\l(|p_n-p_{n+1}\l|+|p_{n+2}-p_{n+3}\l|\l^2+\dots)
\end{eqnarray}
Since $\l\in(1,1+lp_1)\subset(1,1+l\a)$ we have
\begin{equation}\eqlabel{t2eq14}
p_n-p_{n+1}(1+l\a)\leq p_n-p_{n+1}\l\leq p_n-p_{n+1}
\end{equation}
which shows that
\begin{equation}\eqlabel{t2eq15}
|p_n-p_{n+1}\l|\leq p_n-p_{n+1}(1-l\a).
\end{equation}
If we denote by $h=1-l\a$ and we use the bound from Eq.\equref{t2eq15} and the fact that $(p_n)_{n}$ is decreasing in Eq.\equref{t2eq13} we obtain
\begin{eqnarray}\eqlabel{t2eq16}
\frac{|u_n|}{\l}&\leq& p_n-p_{n+1}h+(p_{n+2}-p_{n+3}h)\l^2+\dots\nonumber\\
                &=&(p_n+p_{n+2}\l^2+p_{n+4}\l^4+\dots)-h(p_{n+1}+p_{n+3}\l^2+p_{n+5}\l^4+\dots)\nonumber\\
                &\leq& W-h(p_{n+2}+p_{n+4}\l^2+p_{n+6}\l^4+\dots)\nonumber\\
                &=&W-\frac{h}{\l^2}(W-p_n)=W\left(1-\frac{h}{\l^2}\right)+\frac{h}{\l^2}p_n,
\end{eqnarray}
where by Eq.\equref{t1eq1}
\begin{eqnarray}\eqlabel{t2eq17}
W&=&p_n+p_{n+2}\l^2+p_{n+4}\l^4+\dots\nonumber\\
 &\leq&p_1^{\left[\frac{n+1}{2}\right]}+p_1^{\left[\frac{n+1}{2}\right]}p_1\l^2+p_1^{\left[\frac{n+1}{2}\right]}p_1^2\l^4+\dots\nonumber\\
 &=&p_1^{\left[\frac{n+1}{2}\right]}(1+p_1\l^2+p_1^2\l^4+\dots)=\frac{p_1^{\left[\frac{n+1}{2}\right]}}{1-p_1\l^2}.
\end{eqnarray}
From Eq.\equref{t2eq16} and Eq.\equref{t2eq17} we conclude that
\begin{eqnarray}\eqlabel{t2eq18}
\frac{|u_n|}{\l}&\leq&p_1^{\left[\frac{n+1}{2}\right]}\left[\frac{1}{1-p_1\l^2}+\frac{h}{\l^2}\left(1-\frac{1}{1-p_1\l^2}\right)\right]\nonumber\\
                &=&\frac{1-p_1h}{1-p_1\l^2}p_1^{\left[\frac{n+1}{2}\right]}.
\end{eqnarray}
Until now we have an approximation for $u_n$, but we still need one for $t_n$ and to solve this aspect lets write
\begin{equation}\eqlabel{t2eq19}
T(z)=\frac{1}{U(z)}=\frac{1}{1-(1-U(z))}=\displaystyle\sum_{n\geq0}{(-1)^n(U-1)^n},
\end{equation}
which is true since the convergence of $C(z)$ implies $|z|<\frac{1}{\sqrt{p_1}}$ so that
\begin{eqnarray}\eqlabel{t2eq20}
|1-U|&\leq&\frac{\l(1-p_1h)}{1-p_1\l^2}\displaystyle\sum_{n\geq1}{p_1^{\frac{n}{2}}|z|^n}\leq\frac{|z|\l\sqrt{p_1}\left[1-p_1(1-lp_1)\right]}{(1-p_1\l^2)(1-|z|\sqrt{p_1})}\nonumber\\
     &\leq&\frac{\sqrt{p_1}(1+lp_1^2)(1+lp_1)^2}{\left[1-p_1(1+lp_1)^2\right]\left[1-\sqrt{p_1}(1+lp_1)\right]}<0.8.
\end{eqnarray}
Since $u_0=1$ we have $U-1=\displaystyle\sum_{n\geq1}{u_nz^n}$ and
\begin{equation}\eqlabel{t2eq21}
(U-1)^k=\displaystyle\sum_{l\geq1}{\sum_{\substack{i_1+\dots+i_k=l\\ i_j\geq1,j=\overline{1,k}}}{u_{i_1}\dots u_{i_k}}z^l}.
\end{equation}
Combining Eq.\equref{t2eq19} with Eq.\equref{t2eq21} we get
\begin{equation}\eqlabel{t2eq22}
\displaystyle\sum_{n\geq0}{t_nz^n}=\displaystyle\sum_{k=0}^{\infty}{(-1)^k\sum_{l=1}^{\infty}{b_{l,k}z^l}},
\end{equation}
where
\begin{equation}\eqlabel{t2eq23}
b_{l,k}=\displaystyle\sum_{\substack{i_1+\dots+i_k=l\\ i_j\geq1,j=\overline{1,k}}}{u_{i_1}\dots u_{i_k}}.
\end{equation}
From Eq.\equref{t2eq22} we get that $t_0=1$ and
\begin{equation}\eqlabel{t2eq24}
t_n=\displaystyle\sum_{k=1}^{n}{(-1)^k\sum_{\substack{i_1+\dots+i_k=n\\ i_j\geq1,j=\overline{1,k}}}{u_{i_1}\dots u_{i_k}}},\ \ k\geq1.
\end{equation}
Notice that from Eq.\equref{t2eq18} we can write
\begin{equation}\eqlabel{t2eq25}
|b_{n,k}|\leq \d^k\sum_{\substack{i_1+\dots+i_k=n\\ i_j\geq1,j=\overline{1,k}}}{p_1^{\left[\frac{i_1+1}{2}\right]}\dots p_1^{\left[\frac{i_k+1}{2}\right]}},
\end{equation}
where $\d=\frac{\l(1-p_1h)}{1-p_1\l^2}$. Since, by induction it is easy to verify that
\begin{equation}\eqlabel{t2eq25_2}
\left[\frac{i_1+1}{2}\right]+\dots+\left[\frac{i_k+1}{2}\right]\geq\left[\frac{i_1+\dots+i_k+1}{2}\right]=\left[\frac{n+1}{2}\right],
\end{equation}
and the number of terms in the sum of Eq.\equref{t2eq25} is equal with the number of different positive integers solutions of the equation $i_1+\dots+i_k=n$, which is given by $\binom{n-1}{k-1}$, we deduce
\begin{equation}\eqlabel{t2eq26}
|b_{n,k}|\leq p_1^{\left[\frac{n+1}{2}\right]}\binom{n-1}{k-1}\d^k.
\end{equation}
Now from Eq.\equref{t2eq24} and Eq.\equref{t2eq26} we have
\begin{equation}\eqlabel{t2eq27}
-p_1^{\left[\frac{n+1}{2}\right]}\d\displaystyle\sum_{k=0}^{\left[\frac{n-1}{2}\right]}{\binom{n-1}{2k}\d^{2k}}\leq t_n\leq p_1^{\left[\frac{n+1}{2}\right]}\d\displaystyle\sum_{k=0}^{\left[\frac{n}{2}\right]-1}{\binom{n-1}{2k+1}\d^{2k+1}},
\end{equation}
so that
\begin{equation}\eqlabel{t2eq28}
|t_n|\leq \frac{\d}{2}p_1^{\left[\frac{n+1}{2}\right]}\left[(1+\d)^{n-1}+(1-\d)^{n-1}\right].
\end{equation}
We see from Eq.\equref{t2eq12} and Eq.\equref{t2eq28} that the difference
\begin{eqnarray}\eqlabel{t2eq29}
|q_n\l^n-q_3\l^3|&=&\left|\frac{\displaystyle\sum_{s=0}^{n+1}{t_s\l^s}-\displaystyle\sum_{s=0}^{4}{t_s\l^s}}{\l}\right|=\left|\displaystyle\sum_{s=5}^{n+1}{t_s\l^{s-1}}\right|\nonumber\\
                 &\leq&\frac{\d}{2}\displaystyle\sum_{s=5}^{\infty}{p_1^{\left[\frac{s+1}{2}\right]}\left[(1+\d)^{s-1}+(1-\d)^{s-1}\right]\l^{s-1}}
\end{eqnarray}
If we denote by $\s_1=(1+\d)\l$, $\s_2=(1-\d)\l$ and by $V$ the bound in Eq.\equref{t2eq29}, it is not hard to see that
\begin{equation}\eqlabel{t2eq30}
V=\frac{\d p_1^3}{2}\left[\frac{\s_1^4(1+\s_1)}{1-p_1\s_1^2}+\frac{\s_2^4(1+\s_2)}{1-p_1\s_2^2}\right].
\end{equation}
Recalling that $h=1-lp_1$, $\l\in(1,1+lp_1)$ and $p_1\leq 0.1$ we observe that $\s_2$ is bounded by
\begin{equation*}
-\frac{lp_1[1+2p_1(1+lp_1)]}{1-p_1(1+lp_1)^2}\leq\s_2\leq-\frac{lp_1^2(1+lp_1)}{1-p_1},
\end{equation*}
which gives $|\s_2|<0.5$ and  $\frac{\d\s_2^4(1+\s_2)}{2(1-p_1\s_2^2)}<0.1$.
Substituting the last relation in Eq.\equref{t2eq30} we can rewrite the bound in Eq.\equref{t2eq29} as
\begin{equation}\eqlabel{t2eq31}
|q_n\l^n-q_3\l^3|\leq E(p_1)p_1^3
\end{equation}
where, if we denote by $\et=1+lp_1$,
\begin{equation}\eqlabel{t2eq32}
E(p_1)=0.1+\frac{\et^5\left[1+(1-2p_1)\et\right]^4\left[1+p_1(\et-2)\right]\left[1+\et+(1-3p_1)\et^2\right]}{2(1-p_1\et^2)^4\left[(1-p_1\et^2)^2-p_1\et^2(1+\et-2p_1\et)^2\right]}.
\end{equation}
To obtain $E(\a)$ it is enough to make in Eq.\equref{t2eq32} the substitutions: $p_1\to\a$, $l\to l(\a)$ and $\et\to 1+l\a$. The following lemma gives an approximation for $q_3\l^3$.
\begin{lemma}\lelabel{L1}
If $p_1\leq\a$ then
\begin{equation}\eqlabel{eqL1}
q_3\l^3=1-p_2+2p_3-3p_4+p_1^2+6p_2^2-6p_1p_2+\MO(L(\a)p_1^3)
\end{equation}
where $\MO(x)$ is a function such that $|\MO(x)|\leq|x|$ and $L(\a)$ is an expression depending on $\a$ and given in Eq.\equref{l1eq15}.
\end{lemma}
From Eq.\equref{t2eq29} , Eq.\equref{t2eq32} and \leref{L1} we conclude that
\begin{equation}\eqlabel{t2eq33}
|q_n\l^n-(1-p_2+2p_3-3p_4+p_1^2+6p_2^2-6p_1p_2)|\leq\G(\a)p_1^3,
\end{equation}
with $\G(\a)=L(\a)+E(\a)$ and this ends the proof of the theorem. $\qed$
\subsubsection{Proof of \leref{L1}}
\noindent
From Eq.\equref{eqT1.1} of \thref{T1} we can write
\begin{equation}\eqlabel{l1eq1}
\l=1+p_1-p_2+p_3-p_4+2p_1^2+3p_2^2-5p_1p_2+\MO(K(\a)p_1^3)
\end{equation}
and raising to the third power we get
\begin{eqnarray}\eqlabel{l1eq2}
\l^3&=&(1+\z_1+\z_2)^3+3(1+\z_1+\z_2)^2\MO(K(\a)p_1^3)+\MO(K^3(\a)\a^6p_1^3)+\nonumber\\
    & &+3(1+\z_1+\z_2)\MO(K^2(\a)\a^3p_1^3),
\end{eqnarray}
where we have used the following notations
\begin{eqnarray}
\z_1&=&p_1-p_2,\eqlabel{l1eq3}\\
\z_2&=&p_3-p_4+2p_1^2+3p_2^2-5p_1p_2.\eqlabel{l1eq4}
\end{eqnarray}
Since we can easily see from Eq.\equref{l1eq3} and Eq.\equref{l1eq4} that $\z_1=\MO(p_1)$ and $\z_2=\MO(3p_1^2)$, we deduce that
\begin{equation}\eqlabel{l1eq5}
\l^3-(1+\z_1+\z_2)^3=\MO(S(\a)p_1^3)
\end{equation}
with $S(\a)$ given bellow by
\begin{equation}\eqlabel{l1eq6}
S(\a)=3(1+\a+3\a^2)^2K(\a)+3\a^3(1+\a+3\a^2)K^2(\a)+\a^6K^3(\a).
\end{equation}
Now we observe that by expanding $(1+\z_1+\z_2)^3$ we have
\begin{eqnarray}\eqlabel{l1eq7}
(1+\z_1+\z_2)^3&=&1+3(\z_1+\z_2+\z_1^2)+6\z_1\z_2+3\z_2^2+3\z_1\z_2^2+3\z_1^2\z_2\nonumber\\
               & &+\z_1^3+\z_2^3\nonumber\\
               &=&1+3(\z_1+\z_2+\z_1^2)+\MO(18p_1^3)+\MO(27\a p_1^3)+\MO(27\a^2p_1^3)\nonumber\\
               & &+\MO(9\a p_1^3)+\MO(p_1^3)+\MO(27\a^3p_1^3)
\end{eqnarray}
which along with Eq.\equref{l1eq5} and Eq.\equref{l1eq6} gives the following expression for
\begin{equation}\eqlabel{l1eq8}
\l^3=1+3(\z_1+\z_2+\z_1^2)+\MO(P(\a)p_1^3),
\end{equation}
where $P(\a)$ can be computed by the formula
\begin{eqnarray}\eqlabel{l1eq9}
P(\a)&=&3K(\a)(1+\a+3\a^2)[1+\a+3\a^2+K(\a)\a^3]+\a^6K^3(\a)\nonumber\\
     & &+9\a(4+3\a+3\a^2)+19.
\end{eqnarray}
It is easy to see from Eq.\equref{t2eq5} that
\begin{equation}\eqlabel{l1eq10}
q_3=1-p_1-2(p_1-p_2)+p_1^2-p_3=\MO((1-p_1)^2)
\end{equation}
and combined with Eq.\equref{l1eq8} gives
\begin{equation}\eqlabel{l1eq11}
q_3\l^3=q_3[1+3(\z_1+\z_2+\z_1^2)]+\MO(P(\a)(1-p_1)^2p_1^3).
\end{equation}
The last step in our proof is to find an approximation for the first term on the right in Eq.\equref{l1eq11}. If we write
\begin{equation}\eqlabel{l1eq12}
q_3[1+3(\z_1+\z_2+\z_1^2)]=1-p_2+2p_3-3p_4+p_1^2+6p_2^2-6p_1p_2+H,
\end{equation}
then $H$ checks the relations
\begin{eqnarray}\eqlabel{l1eq13}
H&=&3(p_1-p_2)\left\{(p_1^2-p_3)[3(p_1-p_2)+1-p_2]+p_2^2-9(p_1-p_2)^2\right\}\nonumber\\
 & &-3(p_3-p_4)[p_1+2(p_1-p_2)-(p_1^2-p_3)]\nonumber\\
 &=&\MO(36p_1^3).
\end{eqnarray}
Finally, combining Eq.\equref{l1eq11}, Eq.\equref{l1eq12} and Eq.\equref{l1eq13} we get
\begin{equation}\eqlabel{l1eq14}
q_3\l^3=1-p_2+2p_3-3p_4+p_1^2+6p_2^2-6p_1p_2+\MO(L(\a)p_1^3)
\end{equation}
where we used the notation
\begin{eqnarray}\eqlabel{l1eq15}
L(\a)&=&36+(1-p_1)^2P(\a)\nonumber\\
     &<&3K(\a)(1+\a+3\a^2)[1+\a+3\a^2+K(\a)\a^3]+\a^6K^3(\a)\nonumber\\
     & &+9\a(4+3\a+3\a^2)+55.\ \ \qed
\end{eqnarray}
\subsection{Proof of \coref{C1} and \coref{C2}}
\noindent
To prove the relation in \coref{C1} we see that
\begin{equation}
1+p_1-p_2+2(p_1-p_2)^2=\m-[p_3-p_4-p_2(p_1-p_2)]
\end{equation}
and since
\begin{eqnarray}
p_3-p_4&=&\PP(X_1>x,X_2>x,X_3>x,X_4\leq x)\nonumber\\
        &\leq& p_1\PP(X_1>x,X_2\leq x)=p1(p_1-p_2),
\end{eqnarray}
we get
\begin{equation}
|p_3-p_4-p_2(p_1-p_2)|\leq p_3-p_4+p_2(p_1-p_2)\leq p_1^2.
\end{equation}
Combining these relations with Eq.\equref{eqT1.1} from \thref{T1}, we obtain
\begin{equation}
|\l-(1+p_1-p_2+2(p_1-p_2)^2)|\leq p_1^2(1+K(\a)p_1)\leq (1+\a K(\a)) p_1^2\qed
\end{equation}
\noindent
To prove \coref{C2} notice that
\begin{equation}
-3p_1^2\leq p_1^2+2p_3-3p_4+6p_2^2-6p_1p_2\leq3p_1^2
\end{equation}
and using Eq.\equref{eqT2.1} from \thref{T2}, we get
\begin{eqnarray}
|q_n\l^n-(1-p_2)|&\leq&\G(\a)p_1^3+|p_1^2+2p_3-3p_4+6p_2^2-6p_1p_2|\nonumber\\
                 &\leq&\a\G(\a)p_1^2+3p_1^2=(3+\a\G(\a))p_1^2.\ \ \ \qed
\end{eqnarray}

\subsection{Proof of \thref{T3} and \thref{T4}}
\noindent
Denoting by
\begin{eqnarray*}
\m_1&=&1-p_2+2p_3-3p_4+p_1^2+6p_2^2-6p_1p_2,\\
\m_2&=&1+p_1-p_2+p_3-p_4+2p_1^2+3p_2^2-5p_1p_2
\end{eqnarray*}
we observe that
\begin{equation}
0\leq\frac{\m_1}{\m_2}<1
\end{equation}
and that
\begin{equation}
\m_2=1+\underbrace{(p_1-p_2)(1-p_2)+p_3-p_4+2(p_1-p_2)^2}_{\geq0}\geq1.
\end{equation}
With the help of Eq.\equref{eqT1.1} and Eq.\equref{eqT2.1} we get
\begin{eqnarray}\eqlabel{t3eq1}
\left|q_n-\frac{\m_1}{\m_2^n}\right|&\leq&\left|q_n-\frac{\m_1}{\l^n}\right|+\left|\frac{\m_1}{\l^n}-\frac{\m_1}{\m_2^n}\right|\nonumber\\
                                &\leq&\G(\a)p_1^3+|\m_1|\left|\frac{1}{\l}-\frac{1}{\m_2}\right|\left|\frac{1}{\l^{n-1}}+\dots+\underbrace{\frac{1}{\l^{n-j}\m_2^j}}_{\leq1}+\dots+\frac{1}{\m_2^{n-1}}\right|\nonumber\\
                                &\leq&\left[\G(\a)+nK(\a)\right]p_1^3.
\end{eqnarray}
To express $\m_1$ and $\m_2$ in terms of $q$'s we have to observe first that
\begin{eqnarray*}
    p_1&=&1-q_1\\
    p_2&=&1-2q_1+q_2\\
    p_3&=&1-3q_1+2q_2+q_1^2-q_3\\
    p_4&=&1-4q_1+3q_2-2q_1q_2+3q_1^2-2q_3+q_4
\end{eqnarray*}
and after the proper substitutions, we get
\begin{eqnarray}
\m_1&=&1+q_1-q_2+q_3-q_4+2q_1^2+3q_2^2-5q_1q_2\\
\m_2&=&6(q_1-q_2)^2+4q_3-3q_4
\end{eqnarray}
To finish the proof of \thref{T3} it is enough to use the above relations in Eq.\equref{t3eq1} to obtain
\begin{equation}
\left|q_n-\frac{6(q_1-q_2)^2+4q_3-3q_4}{(1+q_1-q_2+q_3-q_4+2q_1^2+3q_2^2-5q_1q_2)^n}\right|\leq\Delta_1(1-q_1)^3
\end{equation}
where
\begin{equation}
\Delta_1=\Delta_1(\a,n)=\G(\a)+nK(\a). \qed
\end{equation}
For the proof of \thref{T4} we use the same approach as in \thref{T3}. With the help of Eq.\equref{eqC1} and Eq.\equref{eqC2} we get
\begin{eqnarray}\eqlabel{t4eq1}
\left|q_n-\frac{\n_1}{\n_2^n}\right|&\leq&\left|q_n-\frac{\n_1}{\l^n}\right|+\left|\frac{\n_1}{\l^n}-\frac{\n_1}{\n_2^n}\right|\leq(3+p_1\G(\a))p_1^2+n\n_1\frac{|\l-\n_2|}{\l\n_2}\nonumber\\
                                &\leq&\left[3+\G(\a)p_1+n(1+p_1K(\a))\right]p_1^2
\end{eqnarray}
where $\n_1=1-p_2$ and $\n_2=1+p_1-p_2+2(p_1-p_2)^2$. We express $\n_1$ and $\n_2$ in terms of $q$'s using the relations for $p_1$ to $p_4$ from the proof of \thref{T3} and we have
\begin{eqnarray}
\n_1&=&2q_1-q_2\\
\n_2&=&1+q_1-q_2+2(q_1-q_2)^2,
\end{eqnarray}
which we substitute in Eq.\equref{t4eq1} to obtain
\begin{equation}
\left|q_n-\frac{2q_1-q_2}{\left[1+q_1-q_2+2(q_1-q_2)^2\right]^n}\right|\leq\Delta_2(1-q_1)^2
\end{equation}
where
\begin{equation}
\Delta_2=\Delta_2(\a,n,q_1)=3+\G(\a)(1-q_1)+n\left[1+K(\a)(1-q_1)\right]. \qed
\end{equation}


\end{document}